\documentclass{amsart}

\usepackage{fancyhdr}
\usepackage{lastpage}
\usepackage{stmaryrd,yhmath}

\newcommand{\bdis}{\begin{displaymath}}
\newcommand{\edis}{\end{displaymath}}
\newcommand{\be}{\begin{equation}}
\newcommand{\ee}{\end{equation}}

\newcommand{\mcal}{\mathcal}

\newcommand{\vp}{\varphi}

\newcommand{\vth}{\vartheta}

\newcommand{\mG}{\mathring{G}}
\newcommand{\mT}{\mathring{T}}

\newcommand{\zf}{\zeta\left(\frac{1}{2}+it\right)}

\newtheorem{theorem}{Theorem}

\theoremstyle{definition}

\theoremstyle{remark}
\newtheorem{remark}[]{Remark}

\numberwithin{equation}{section}

\begin{document}

\title{Jacob's ladders, heterogeneous quadrature formulae, big asymmetry and related formulae for the Riemann zeta-function}

\author{Jan Moser}

\address{Department of Mathematical Analysis and Numerical Mathematics, Comenius University, Mlynska Dolina M105, 842 48 Bratislava, SLOVAKIA}

\email{jan.mozer@fmph.uniba.sk}


\begin{abstract}
In this paper we obtain as our main result new class of formulae expressing correlation integrals of the third-order in $Z$ on disconnected
sets $\mathring{G}_1(x),\mathring{G}_2(y)$ by means of an autocorrelative sum of the second order in $Z$. Moreover, the distance of the sets
$\mathring{G}_1(x),\mathring{G}_2(y)$ from the set of arguments of autocorrelative sum is extremely big, namely $\sim A\pi(T),\ T\to\infty$, where
$\pi(T)$ is the prime-counting function.
\end{abstract}
\maketitle

\section{Introduction}

\subsection{}

Let us remind the following formulae
\be \label{1.1}
\begin{split}
 & \int_T^{2T} Z^4(t){\rm d}t\sim \frac{2\pi}{\ln T}\sum_{T\leq t_\nu\leq T+2T} Z^4(t_\nu),\ T\to\infty , \\
 & \int_T^{T+U} Z^2(t){\rm d}t\sim \frac{2\pi}{\ln T}\sum_{T\leq t_\nu\leq T+U} Z^2(t_\nu),\ U=\sqrt{T}\ln T,\ T\to\infty
\end{split}
\ee
(see \cite{5}, (4.4), (4.7)), where
\be \label{1.2}
\begin{split}
 & Z(t)=e^{i\vth(t)}\zf, \\
 & \vth(t)=-\frac t2\ln \pi +\text{Im}\ln\Gamma\left(\frac 14+i\frac t2\right)
\end{split}
\ee
(see \cite{9}, pp. 79, 329), and $\{ t_\nu\}$ is the Gram sequence (comp. \cite{9}, p. 99). The formulae (\ref{1.1}) were proved by us in connection with the
Kotelnikov-Whittaker-Nyquist theorem. Namely, by these formulae we have expressed the biquadratic and quadratic effects for the continuous signals
\bdis
Z(t),\quad t\in [T,2T];\ t\in [T,T+U]
\edis
from the point of view of information theory.

\subsection{}

The formulae (\ref{1.1}) are:
\begin{itemize}
 \item[(a)] asymptotic quadrature formulae (from the left to the right),
 \item[(b)] asymptotic summation formulae (from the right to the left).
\end{itemize}

Next, for the second formula in (\ref{1.1}) (for example) we have
\be \label{1.3}
\begin{split}
 & t\in [T,T+U] \rightarrow t_\nu\in [T,T+U], \\
 & Z^2(t)\rightarrow Z^2(t_\nu),
\end{split}
\ee
i. e. the segment $[T,T+U]$ and the exponent $2$ are conserved.

\begin{remark}
By (\ref{1.3}) it is natural to call the formulae of the kind (\ref{1.1}) \emph{homogeneous formulae}.
\end{remark}

On the contrary, we obtain in this paper some heterogeneous asymptotic quadrature formulae -- such formulae that the properties (\ref{1.3}) are not
fulfilled.

\section{Heteregeneous quadrature formulae}

\subsection{}

Let (see \cite{3}, (3))
\be \label{2.1}
\begin{split}
 & G_1(x)=G_1(x;T,H)= \\
 & = \bigcup_{T\leq t_{2\nu}\leq T+H}\{ t:\ t_{2\nu}(-x)<t<t_{2\nu}(x)\},\quad 0<x\leq \frac{\pi}{2}, \\
 & G_2(y)=G_2(y;T,H)= \\
 & = \bigcup_{T\leq t_{2\nu+1}\leq T+H}\{ t:\ t_{2\nu+1}(-y)<t<t_{2\nu+1}(y)\}, \quad 0<y\leq \frac{\pi}{2}, \\
 & H=T^{1/6+2\epsilon},
\end{split}
\ee
where the collection of sequences
\bdis
\{ t_\nu(\tau)\},\quad \tau\in [-\pi,\pi],\ \nu=1,2,\dots
\edis
is defined by the equation (see \cite{3}, (1))
\bdis
\vth[t_\nu(\tau)]=\pi\nu+\tau;\ t_\nu(0)=t_\nu,
\edis
and (see \cite{3}, (8))
\be \label{3.2}
m\{ G_1(x)\}=\frac{x}{\pi}H+\mcal{O}\left(\frac{x}{\ln T}\right),\
m\{ G_2(y)\}=\frac{y}{\pi}H+\mcal{O}\left(\frac{y}{\ln T}\right),
\ee
where $m\{ G_1(x)\}$, \dots is the measure of the set $G_1$, \dots. \\

Let
\bdis
\vp_1\{ \mG_1(x)\}=G_1(x),\quad \vp_1\{\mG_2(y)\}=G_2(y).
\edis
The following Theorem holds true.

\begin{theorem}
 \be \label{2.3}
 \begin{split}
  & \int_{\mG_1(x)}\omega(t)Z[\vp_1(t)]Z^2(t){\rm d}t = \\
  & = \frac{m\{ G_1(x)\}}{Q_1\ln P_0}\sum_{T\leq t_\nu\leq T+U_1}Z(t_\nu)Z\left( t_\nu+\frac{x}{\ln P_0}\right)+\mcal{O}\left(\frac{H}{\ln T}\right), \\
  & \int_{\mG_2(y)}\omega(t)Z[\vp_1(t)]Z^2(t){\rm d}t = \\
  & = -\frac{m\{ G_2(y)\}}{Q_1\ln P_0}\sum_{T\leq t_\nu\leq T+U_1}Z(t_\nu)Z\left( t_\nu+\frac{y}{\ln P_0}\right)+\mcal{O}\left(\frac{H}{\ln T}\right), \\
  & x,y\in \left(\left. 0,\frac{\pi}{2}\right.\right],\quad T\to\infty,
 \end{split}
\ee
where
\be \label{2.4}
\omega(t)=\frac{1}{\ln t}\left\{ 1+\mcal{O}\left(\frac{\ln\ln t}{\ln t}\right)\right\},
\ee
and (see \cite{1}, (38); $H\to U_1$)
\be \label{2.5}
\begin{split}
 & Q_1=\sum_{T\leq t_\nu\leq T+U_1}1=\frac 1\pi U_1\ln P_0+\mcal{O}\left(\frac{U_1^2}{T}\right), \\
 & U_1=\sqrt{T}\ln P_0,\ P_0=\sqrt{\frac{T}{2\pi}}.
\end{split}
\ee
Furthermore we have
\be \label{2.6}
m\{\mG_1(x)\},\ m\{\mG_2(y)\}<T^{1/3+\epsilon}.
\ee
\end{theorem}

\subsection{}

\begin{remark}
 Every of the correlation integrals in (\ref{2.3}) contains the product
 \be \label{2.7}
 Z[\vp_1(t)]Z^2(t)
 \ee
and, as usually, we have (see \cite{6}, (6.2); \cite{7}, (8.4))
\be \label{2.8}
t-\vp_1(t)\sim (1-c)\frac{t}{\ln t}\sim (1-c)\pi(t),\quad t\to\infty,
\ee
i. e. we have big difference of arguments in (\ref{2.7}), where $c$ is the Euler constant and $\pi(t)$ is the prime-counting function.
\end{remark}

Next, in the case
\bdis
\begin{split}
 & t=\mT;\ \vp_1(\mT)=T \ \Rightarrow \ T\to\infty \ \Leftrightarrow \ \mT\to\infty
\end{split}
\edis
we obtain from (\ref{2.8})
\bdis
\begin{split}
 & \mT-T\sim (1-c)\frac{\mT}{\ln\mT} \ \Rightarrow \ 1-\frac{T}{\mT}\sim (1-c)\frac{1}{\ln \mT} \ \Rightarrow \\
 & \Rightarrow\ \mT \sim T \ \Rightarrow \ \ln \mT\sim \ln T,
\end{split}
\edis
i. e.
\bdis
\mT-T\sim (1-c)\frac{T}{\ln T},
\edis
and (see (\ref{2.5}))
\bdis
\mT-(T+U_1)\sim (1-c)\frac{T}{\ln T}-U_1\sim (1-c)\frac{T}{\ln T}.
\edis
Consequently we have
\be \label{2.9}
\begin{split}
 & \rho\{ [\mT,\widering{T+H}];[T,T+U_1]\}\sim (1-c)\pi(T); \\
 & [T,T+U_1]\prec [\mT,\widering{T+H}],
\end{split}
\ee
where $\rho$ stands for the distance of the corresponding segments. We may, of course, put
\be \label{2.10}
G_1(x)\cap [T,T+H]=\bar{G}_1(x)\rightarrow G_1(x),\dots
\ee
if necessary.

\begin{remark}
We have the following properties
\begin{itemize}
 \item[(a)]
 \bdis
 \mG_1(x),\mG_2(y)\in [\mT,\widering{T+H}]\rightarrow [T,T+U],
 \edis
 (comp. (\ref{1.3}), (\ref{2.10})), where $\mG_1(x),\mG_2(y)$ are disconnected sets,
 \item[(b)] if
 $$\mcal{G}=\mG_1(x),\mG_2(y)$$
 then extremely big distance occurs, namely (comp. (\ref{2.9}), (\ref{2.10}))
 \bdis
 \rho\{\mcal{G};[T,T+U_1]\}\sim (1-c)\pi(T),\ T\to\infty,
 \edis
 \item[(c)] for the corresponding orders of $Z$ (comp. (\ref{1.3}))
 \bdis
 1+2\rightarrow 1+1
 \edis
 then the formulae (\ref{2.3}) are \emph{strongly heterogeneous} (comp. Remark 1).
\end{itemize}

\end{remark}

\begin{remark}
Moreover we explicitly notice the following:
\begin{itemize}
 \item[(a)] the formulae (\ref{2.3}) are not accessible by the current methods in the theory of the Riemann zeta-function,
 \item[(b)] small improvements of the exponents
 \bdis
 \frac 16,\ \frac 12, \ \dots
 \edis
 are irrelevant for main direction of this paper (comp. our paper \cite{2}, Appendex A: On I.M. Vinogradov' scepticism on possibilities of the method
 of trigonometric sums).
\end{itemize}
\end{remark}

\section{Big asymmetry and related formulae}

\subsection{}

Let (see \cite{4}, p. 29)
\be \label{3.1}
\begin{split}
 & G_3(x)=G_3(x;T,U_2)= \\
 & = \bigcup_{T\leq g_{2\nu}\leq T+U_2}\{ t:\ g_{2\nu}(-x)<t<g_{2\nu}(x)\},\quad 0<x\leq \frac{\pi}{2}, \\
 & G_4(y)=G_4(y;T,U_2)= \\
 & = \bigcup_{T\leq g_{2\nu+1}\leq T+U_2}\{ t:\ g_{2\nu+1}(-y)<t<g_{2\nu+1}(y)\}, \quad 0<y\leq \frac{\pi}{2}, \\
 & U_2=T^{5/12+2\epsilon},
\end{split}
\ee
where the collection of sequences
\bdis
\{ g_\nu(\tau)\},\ \tau\in [-\pi,\pi],\ \nu=1,2,\dots
\edis
is defined by the equation (see \cite{4}, (6))
\bdis
\vth_1[g_\nu(\tau)]=\frac{\pi}{2}\nu+\frac{\tau}{2};\ g_\nu(0)=g_\nu,
\edis
where (comp. (\ref{1.2}))
\bdis
\vth(t)=\vth_1(t)+\mcal{O}\left(\frac 1t\right),\ \vth_1(t)=\frac t2\ln\frac{t}{2\pi}-\frac t2-\frac{\pi}{8},
\edis
and (see \cite{4}, (13))
\be \label{3.3}
\begin{split}
 & m\{ G_3\}=\frac{x}{\pi}U_2+\mcal{O}\left(\frac{x}{\ln T}\right), \\
 & m\{ G_4\}=\frac{y}{\pi}U_2+\mcal{O}\left(\frac{y}{\ln T}\right).
\end{split}
\ee
Let
\bdis
\vp_1\{\mG_3(x)\}=G_3(x),\quad \vp_1\{\mG_4(y)\}=G_4(y).
\edis
The following Theorem holds true.

\begin{theorem}
\be \label{3.4}
\begin{split}
 & \int_{\mG_3(x)}\omega(t)Z^2[\vp_1(t)]Z^2(t){\rm d}t\sim \\
 & \sim \frac{4}{\pi}U_2\sin x+\int_{\mG_4(x)}\omega(t)Z^2[\vp_1(t)]Z^2(t){\rm d}t, \\
 & x\in \left(\left. 0,\frac{\pi}{2}\right.\right],\quad T\to\infty.
\end{split}
\ee
\end{theorem}

\begin{remark}
By the asymptotic formula (\ref{3.4}) is expressed the property of big asymmetry in the distribution of the values of $Z$ on disconnected sets
$\mG_3(x),\mG_4(x)$. Namely, the correlation integral (comp. (\ref{2.8})) on the disconnected set $\mG_3(x)$ essentially exceeds of that on the set
$\mG_4(x)$. For example
\bdis
\int_{\mG_3(\pi/2)}-\int_{\mG_4(\pi/2)}\sim \frac{4}{\pi}U_2,\ T\to\infty,
\edis
where
\bdis
m\{\mG_3(\pi/2)\}+m\{\mG_4(\pi/2)\}=U_2
\edis
if (comp. (\ref{2.10}))
\bdis
G_3(x)\cap [T,T+U_2] \rightarrow G_3(x), \dots
\edis
\end{remark}

\subsection{}

Finally, the following Theorem holds true.

\begin{theorem}
\be \label{3.5}
\begin{split}
 & \int_{\mG_3(x)}\omega(t)Z[\vp_1(t)]Z^2(t){\rm d}t= \\
 & =\frac{m\{ G_1(x)\}}{2m\{ G_3(x)\}}\int_{\mG_3(x)}\omega(t)Z^2[\vp_1(t)]Z^2(t){\rm d}t- \\
 & - \frac{m\{ G_1(x)\}}{2m\{ G_4(x)\}}\int_{\mG_4(x)}\omega(t)Z^2[\vp_1(t)]Z^2(t){\rm d}t+\mcal{O}(HT^{-\epsilon}), \\
 & \int_{\mG_2(y)}\omega(t)Z[\vp_1(t)]Z^2(t){\rm d}t= \\
 & =\frac{m\{ G_2(y)\}}{2m\{ G_4(y)\}}\int_{\mG_4(y)}\omega(t)Z^2[\vp_1(t)]Z^2(t){\rm d}t- \\
 & - \frac{m\{ G_2(y)\}}{2m\{ G_3(y)\}}\int_{\mG_3(y)}\omega(t)Z^2[\vp_1(t)]Z^2(t){\rm d}t+\mcal{O}(HT^{-\epsilon}), \\
 & x,y\in \left(\left. 0,\frac{\pi}{2}\right.\right],\ T\to\infty .
\end{split}
\ee
\end{theorem}

\begin{remark}
For the disconnected sets
\bdis
\mG_1(x),\mG_2(y),\mG_3(x),\mG_4(x),\mG_3(y),\mG_4(y);\mG_1\cap\mG_2=\emptyset, \mG_3\cap\mG_4=\emptyset
\edis
we have the following property: the correlation integrals of the order $1+2$ on $\mG_1(x),\mG_2(y)$ are expressed as the linear combinations of the
correlation integrals of the order $2+2$ on $\mG_3(x),\mG_4(x);\mG_3(y),\mG_4(y)$ correspondingly.
\end{remark}

\section{Proof of Theorem 1}

\subsection{}

In the paper \cite{1} (see (10)) we have proved the following autocorrelative formula
\be \label{4.1}
\begin{split}
 & \sum_{T\leq t_\nu\leq T+U_1}Z(t_\nu)Z(t_\nu+\beta)= \\
 & = \frac{2}{\pi}\frac{\sin(\beta\ln P_0)}{\beta\ln P_0}U_1\ln P_0+\mcal{O}(\sqrt{T}\ln^2T),
\end{split}
\ee
where
\be \label{4.2}
\beta=\mcal{O}\left(\frac{1}{\ln T}\right),\quad U_1=\sqrt{T}\ln P_0.
\ee
In the case
\bdis
\beta\ln P_0=x \ \Rightarrow \ \beta=\frac{x}{\ln P_0},\ x\in \left(\left. 0,\frac{\pi}{2}\right.\right]
\edis
we obtain from (\ref{4.1}), (\ref{4.2})
\be \label{4.3}
\begin{split}
 & \sum_{T\leq t_\nu\leq T+U_1}Z(t_\nu)Z\left( t_\nu+\frac{x}{\ln P_0}\right)= \\
 & = \frac{2}{\pi}\frac{\sin x}{x}U_1\ln^2P_0+\mcal{O}(\sqrt{T}\ln^2T).
\end{split}
\ee
Consequently we have (see (\ref{2.5}), (\ref{4.2}))
\be \label{4.4}
\begin{split}
 & \frac{1}{Q_1\ln P_0}\sum_{T\leq t_\nu\leq T+U_1}Z(t_\nu)Z\left( t_\nu+\frac{x}{\ln P_0}\right)= \\
 & = 2\frac{\sin x}{x}+\mcal{O}\left(\frac{1}{\ln T}\right);\quad x\in \left(\left. 0,\frac{\pi}{2}\right.\right] \ \Rightarrow \
 \frac{\sin x}{x}\in \left[\left. \frac{2}{\pi},1\right)\right. .
\end{split}
\ee

\subsection{}

Next, in the paper \cite{3}, (5),(9) we have obtained the following mean-value formula on the disconnected sets $G_1(x),G_2(y)$
(see (\ref{2.1}))
\be \label{4.5}
\begin{split}
 & \int_{G_1(x)}Z(t){\rm d}t=\frac{2}{\pi}H\sin x+\mcal{O}(xT^{1/6+\epsilon}), \\
 & \int_{G_2(y)}Z(t){\rm d}t=-\frac{2}{\pi}H\sin y+\mcal{O}(yT^{1/6+\epsilon}), \\
 & H=T^{1/6+2\epsilon};\ T^{1/6}\psi^2\ln^5T\to T^{1/6+\epsilon}.
\end{split}
\ee
Hence, from (\ref{4.4}) by (\ref{2.3}) we obtain
\be \label{4.6}
\begin{split}
 & \frac{1}{m\{ G_1(x)\}}\int_{G_1(x)}Z(t){\rm d}t=2\frac{\sin x}{x}+\mcal{O}(T^{-\epsilon}), \\
 & \frac{1}{m\{ G_2(y)\}}\int_{G_2(y)}Z(t){\rm d}t=-2\frac{\sin y}{y}+\mcal{O}(T^{-\epsilon}).
\end{split}
\ee

\subsection{}

In the paper \cite{7}, (9.2), (9.5) we have proved the following Lemma: if
\bdis
\vp_1\{ [\mT,\widering{T+U}]\}=[T,T+U],
\edis
then for every Lebesgue-integrable function
\bdis
f(x),\ x\in [T,T+U]
\edis
we have
\be \label{4.7}
\begin{split}
 & \int_{\mT}^{\widering{T+U}} f[\vp_1(t)]\tilde{Z}^2(t){\rm d}t=\int_T^{T+U} f(x){\rm d}x, \\
 & T\geq T_0[\vp_1],\ U\in \left(\left. 0,\frac{T}{\ln T}\right]\right. ,
\end{split}
\ee
where
\be \label{4.8}
\begin{split}
 & \tilde{Z}^2(t)=\frac{Z^2(t)}{\left\{ 1+\mcal{O}\left(\frac{\ln\ln t}{\ln t}\right)\right\}\ln t}=
 \omega(t)Z^2(t); \\
 & \omega(t)=\frac{1}{\left\{ 1+\mcal{O}\left(\frac{\ln\ln t}{\ln t}\right)\right\}\ln t}=
 \frac{1}{\ln t}\left\{ 1+\mcal{O}\left(\frac{\ln\ln t}{\ln t}\right)\right\}.
\end{split}
\ee
Consequently, we have (see (\ref{4.7}), (\ref{4.8}))
\be \label{4.9}
\begin{split}
 & \int_{\mT}^{\widering{T+U}} \omega(t)f[\vp_1(t)]Z^2(t){\rm d}t=\int_T^{T+U} f(x){\rm d}x, \\
 & U\in \left(\left. 0,\frac{T}{\ln T}\right]\right. .
\end{split}
\ee

\subsection{}

Hence, by (\ref{4.6}), (\ref{4.9}) we obtain following formulae
\be \label{4.10}
\begin{split}
 & \frac{1}{m\{ G_1(x)\}}\int_{\mG_1(x)}\omega(t)Z[\vp_1(t)]Z^2(t){\rm d}t= \\
 & = 2\frac{\sin x}{x}+\mcal{O}(T^{-\epsilon}), \\
 & \frac{1}{m\{ G_2(y)\}}\int_{\mG_2(y)}\omega(t)Z[\vp_1(t)]Z^2(t){\rm d}t= \\
 & =-2\frac{\sin y}{y}+\mcal{O}(T^{-\epsilon}) .
\end{split}
\ee
Finally, simple elimination of the values
\bdis
2\frac{\sin x}{x},\ 2\frac{\sin y}{y}
\edis
from (\ref{4.4}), (\ref{4.9}) gives (\ref{2.3}).

\subsection{}

In the case
\bdis
f(x)=1
\edis
in (\ref{4.7}) we obtain (see (\ref{2.1}), (\ref{2.10}), (\ref{4.8}))
\bdis
\begin{split}
 & \int_{\mG_1(x)}\left|\zf\right|^2{\rm d}t<\int_{\mT}^{\widering{T+H}}\left|\zf\right|^2{\rm d}t\sim \\
 & \sim \ln T\int_T^{T+H} 1\cdot {\rm d}t,
\end{split}
\edis
i. e.
\be \label{4.11}
\begin{split}
 & \int_{\mT}^{\widering{T+H}}\left|\zf\right|^2{\rm d}t< A H\ln T,\ T\to\infty .
\end{split}
\ee
Let
\be \label{4.12}
\widering{T+H}-\mT\geq T^{1/3+\epsilon},
\ee
(comp. \cite{6}, (2.5); where $\frac 13$ is the Balasubramanian exponent). Then
\be \label{4.13}
\int_{\mT}^{\widering{T+H}}\left|\zf\right|^2{\rm d}t\sim (\widering{T+H}-\mT)\ln T> B (\widering{T+H}-\mT)\ln T,
\ee
and (see (\ref{4.11}), (\ref{4.13}))
\be \label{4.14}
\widering{T+H}-\mT< CH=CT^{1/6+\epsilon}.
\ee
Now, (\ref{4.14}) contradicts (\ref{4.12}). Consequently we have that
\bdis
m\{\mG_1(x)\}<\widering{T+H}-\mT<T^{1/3+\epsilon}
\edis
and we obtain the second inequality in (\ref{2.6}) by the similar way.

\section{Proofs of Theorem 2 and Theorem 3}

\subsection{}

In the paper \cite{4}, (14), (15) we have proved the following mean-value formulae
\be \label{5.1}
\begin{split}
 & \int_{G_3(x)}Z^2(t){\rm d}t= \\
 & = \frac x\pi U_2\ln\frac{T}{2\pi}+\frac{2}{\pi}( cx+\sin x)U_2+\mcal{O}(T^{5/12}\ln^2T), \\
 & \int_{G_4(y)}Z^2(t){\rm d}t= \\
 & = \frac y\pi U_2\ln\frac{T}{2\pi}+\frac{2}{\pi}( cy-\sin y)U_2+\mcal{O}(T^{5/12}\ln^2T)
\end{split}
\ee
on the disconnected sets $G_3(x),G_4(y)$, (comp. (\ref{3.1})). From (\ref{5.1}) we have in the case $x=y$
\bdis
\begin{split}
 & \int_{G_3(x)} Z^2(t){\rm d}t-\int_{G_4(x)}Z^2(t){\rm d}t= \\
 & =\frac{4}{\pi}U_2\sin x+\mcal{O}(x T^{5/12}\ln^2T).
\end{split}
\edis
Consequently, we obtain from (\ref{5.1}) by (\ref{4.7}) -- (\ref{4.9}) with $f(x)=Z^2(x)$ the formula (\ref{3.4}).

\subsection{}

Next, from the formula (see \cite{4}, (16))
\bdis
\begin{split}
 & \frac{1}{m\{ G_3(x)\}}\int_{G_3(x)}Z^2(t){\rm d}t-\frac{1}{m\{ G_4(x)\}}\int_{G_4(x)}Z^2(t){\rm d}t\sim \\
 & \sim 4\frac{\sin x}{x},\ T\to\infty
\end{split}
\edis
we obtain by (\ref{4.7}) -- (\ref{4.9}) ($f(x)=Z^2(x)$) the following formula
\be \label{5.2}
\begin{split}
 & \frac{1}{m\{ G_3(x)\}}\int_{\mG_3(x)}\omega(t)Z^2[\vp_1(t)]Z^2(t){\rm d}t-\\
 & - \frac{1}{m\{ G_4(x)\}}\int_{\mG_4(x)}\omega(t)Z^2[\vp_1(t)]Z^2(t){\rm d}t\sim 4\frac{\sin x}{x},\ T\to\infty
\end{split}
\ee
and, by the similar way, we obtain the formula for $y\in (0,\pi/2]$. Hence, the elimination of
\bdis
2\frac{\sin x}{x},\ 2\frac{\sin y}{y}
\edis
from (\ref{4.10}), (\ref{5.2}) implies the formula (\ref{3.5}).

\thanks{I would like to thank Michal Demetrian for his help with electronic version of this paper.}

\end{document}